\documentclass[11pt,twoside]{amsart}
\usepackage[utf8]{inputenc}
\usepackage[a4paper,margin=3cm]{geometry}
\usepackage[pdftex,pdfborder={0 0 0},colorlinks,pagebackref]{hyperref} 
\usepackage{latexsym,lmodern,graphicx,enumitem,amssymb,amsfonts,amsmath,amsthm,dsfont,color}

\newtheorem{theorem}{Theorem}
\newtheorem{proposition}[theorem]{Proposition}

\newtheorem{corollary}[theorem]{Corollary}

\theoremstyle{definition}
\newtheorem{example}[theorem]{Example}

\newtheorem{remark}[theorem]{Remark}

\begin{document}

\title{Bakry-Emery calculus for diffusion with additional multiplicative term}
\author{C. Roberto, B. Zegarlinski}

\thanks{Supported by the grants ANR-15-CE40-0020-03 - LSD - Large Stochastic Dynamics, ANR 11-LBX-0023-01 - Labex MME-DII 
and Fondation Simone et Cino del Luca in France, 
and the grant ... in the UK.}

\address{Universit\'e Paris Nanterre, Modal'X, FP2M, CNRS FR 2036, 200 avenue de la R\'epublique 92000 Nanterre, France}
\address{Imperial College of London, Faculty of Natural Sciences, Department of Mathematics,
Huxley Building, South Kensington Campus, London SW7 2AZ, UK}

\email{croberto@math.cnrs.fr, b.zegarlinski@imperial.ac.uk}

\date{\today}

\begin{abstract}
We extend the $\Gamma_2$ calculus of Bakry and Emery to include a \textit{Carr\'e du champ} operator with multiplicative term, providing   results which allow to analyse inhomogeneous diffusions. 
\end{abstract}

\maketitle


The aim of this paper is to extend Bakry-Emery approach \cite{bakry-emery,bakry} to deal with some quantities involving operators not only of order one, but also including order zero. One of the motivations for that is a possible application to analysis of hypercontractivity properties for some classes of inhomogeneous Markov semi-groups, see e.g. \cite{RZ}.

The setting is as follows: $(Q_t)_{t \geq 0}$ is the semi-group associated to a diffusion operator
$L=\Delta - \nabla U \cdot \nabla$, on $\mathbb{R}^n$, where the dot sign stands for the Euclidean scalar product.   We assume that $U \colon \mathbb{R}^n \to \mathbb{R}$   
is twice differentiable
and satisfies
$\int e^{-U} dx = 1$ so that   $\mu(dx)=e^{-U(x)}dx$ is a probability measure on $\mathbb{R}^n$.
By construction $L$ is symmetric in $\mathbb{L}_2(\mu)$. 
Following Bakry-Emery, we denote by $\Gamma$ the \emph{carr\'e du champs} bilinear form
$$
\Gamma(f,g):= \frac{1}{2}\left( L(fg) - f Lg - gLf \right)
$$
and set $\Gamma(f):=\Gamma(f,f)$. For the diffusion $L$ considered here, we have 
$\Gamma(f,g)=\nabla f \cdot \nabla g$. The iterated operator $\Gamma_2$ is defined as
$$
\Gamma_2(f,g)=\frac{1}{2} \left( L\Gamma(f,g)- \Gamma(Lf,g) - \Gamma(f,Lg) \right),
$$ 
and again, for simplicity, we set $\Gamma_2(f):=\Gamma(f,f)$. 
One can see that $\mathrm{Hess}(U) \geq \rho$ (as a matrix), $\rho \in \mathbb{R}$, implies
$\Gamma_2(f) \geq \rho \Gamma(f)$ for all smooth enough $f$ (see \textit{e.g.}\ \cite[Chapter 5]{ane}, \cite{BGL}).

One fundamental result of Bakry and Emery \cite{bakry-emery} is that $\Gamma_2 \geq \rho \Gamma$ (the so-called $\Gamma_2$-condition) is equivalent to the following commutation property between the semi-group and the gradient operator (equivalently $\Gamma$): 
\begin{equation} \label{eq:be}
\Gamma(Q_tf) \leq e^{-2\rho t} Q_t(\Gamma(f)), \qquad t \geq 0
\end{equation}
for any $f$ for which $\Gamma(f)$ is well defined, see \textit{e.g.}\ \cite[Proposition 5.4.1]{ane}. 

We will prove that a similar equivalence holds for an extended operator we introduce now.
Let $W \colon \mathbb{R}^n \to \mathbb{R}_+$ be smooth enough so that in particular $W^2$ 
 belongs to the domain of $L$. Then, for $f,g$ smooth enough, we set
$$
\Gamma^W(f,g):=\Gamma(f,g)+W^2fg
$$
which is therefore a positive bilinear form. The operator $\Gamma^W$ acts as a derivative and multiplicatively. In fact, one can see $\Gamma^W$ as  two-dimensional operator, call it $D$, that acts as $Df:=(\nabla f, Wf)$. With this notation,
$\Gamma^W(f)=|Df|^2=|\nabla f|^2 + W^2f^2$ is nothing but the Euclidean norm squared of the 2 dimensional vector $Df$.
 Similarly, we introduce the iterated operator
\begin{align*}
\Gamma_2^W(f,g) 
& 
:= \frac{1}{2} \left( L\Gamma^W(f,g)- \Gamma^W(Lf,g) - \Gamma^W(f,Lg) \right) \\
& 
= \Gamma_2(f,g) + \frac{1}{2}fgL(W^2) + W^2 \Gamma(f,g) + 2 W \nabla W \nabla (fg) 
\end{align*}
 where the last equality follows from some algebra.

 It should be clear from the definition that we are not dealing with $\Gamma_2$ calculus for the operator $L^W:=L-2W^2$, even though by simple algebra we have $\frac{1}{2}\left( L^W(fg) - f L^Wg - gL^Wf  \right) = \Gamma^W(f,g)$. The point is that we want to derive commutation formulas for the semi-group
$(Q_t)_{t \geq 0}$ associated to $L$, and not for the semi-group associated to $L^W$. This also explains why $\Gamma_2^W$ is defined through the operator $L$ and not $L^W$.\\

\noindent Our first main  result reads as follows. 

\begin{theorem} \label{th:g2w}
Let $\rho \in \mathbb{R}$. The following are equivalent\\

$(i)$ for all $f$ smooth enough $\Gamma_2^W(f) \geq \rho \Gamma^W(f)$\\

$(ii)$ for all $f$ smooth enough and all $t \geq 0$,
$$
\Gamma^W(Q_t f) \leq e^{-2\rho t} Q_t(\Gamma^W(f))
$$
$(iii)$ for all $f$ smooth enough and all $t \geq 0$, 
$$
Q_t(f^2) - (Q_t f)^2 + 2 \int_0^t Q_s \left( W^2 (Q_{t-s} f)^2 \right) ds \leq \frac{1-e^{-2\rho t}}{\rho} Q_t (\Gamma^W(f)) .
$$
\end{theorem}

\begin{remark}
Above, when $\rho=0$, the ratio $\frac{1-e^{-2\rho t}}{\rho}$ is understood as its limit (i.e. $2t$). Notice that it is always non-negative.

Observe that, applying $(ii)$ to constant functions $f \equiv C$, $C \neq 0$, leads to $W^2 \leq e^{-2\rho t} Q_t(W^2)$. Therefore, if $\int W^2 d\mu < \infty$
and $\rho >0$,
 taking the limit $t \to \infty$ and by ergodicity, we would conclude that $W \equiv 0$. Therefore, for the inequality 
 $\Gamma_2^W(f) \geq \rho \Gamma^W(f)$ to hold for a non  trivial $W$, either $\rho \leq 0$ or $\int W^2 d\mu = \infty$. But we have no this restriction removing mean value $\mu(f)$ of the function $f$.

Taking the mean with respect to $\mu$ in $(iii)$ and passing to the limit $t \to \infty$, we get by invariance and ergodicity that
for $\rho >0$, it holds
$$
\int f^2 d\mu - \left( \int f d\mu \right)^2 - 2 \int W^2 \int_0^\infty (Q_s f)^2ds d\mu \leq \frac{1}{\rho} \left( \int |\nabla f|^2 d\mu + \int f^2 W^2 d\mu \right) .
$$
This is a sort of Poincar\'e inequality in particular for a function $f$ with mean value zero. Following Bakry-Emery (see \textit{e.g.}\ \cite[proposition 5.5.4]{ane}), one can actually proove that the latter holds under the weaker assumption that $\int \Gamma_2^W(f) d\mu \geq \rho \int \Gamma^W(f) d\mu$.
\end{remark}

\begin{proof}
The proof mimics the usual case ($W=0$).

To prove that $(i)$ implies $(ii)$, fix $t >0$ and consider the following function
$\Psi \colon s \in [0,t] \to \Psi(s) = Q_s \left( \Gamma^W(Q_{t-s}f) \right) = Q_s (\Gamma(Q_{t-s}f) +W^2 (Q_{t-s}f)^2)$. Then, setting $g:=Q_{t-s}f$, it holds
\begin{align*}
\Psi'(s) 
& = Q_s \left( L \Gamma(g) + L(W^2g^2) - 2 \Gamma(g,Lg)
- 2W^2 g Lg \right) \\
&=
2 Q_s \left( \Gamma_2^W(g) \right) \geq 2\rho Q_s \left( \Gamma^W(g) \right) = 2\rho \Psi(s)
\end{align*}
from which the result of Item $(ii)$ follows.

Now we prove that $(ii)$ implies $(iii)$. Let 
$$
\Psi(s):=Q_s((Q_{t-s}f)^2) +  2\int_0^s Q_u \left( W^2 (Q_{t-u} f)^2 \right) du, \qquad s \in [0,t]. 
$$
Then, setting again $g:=Q_{t-s}f$, it holds
\begin{align*}
\Psi'(s)
=
Q_s \left( L(g^2) -2gLg+ 2W^2(Q_{t-s}f)^2 \right) 
= 
2 Q_s ( \Gamma^W(g) ).
\end{align*}
Therefore
\begin{align*}
\Psi(t) - \Psi(0) 
& = 
\int_0^t \Psi'(s)ds 
=
2 \int_0^t Q_s \left( \Gamma^W(Q_{t-s}f) \right) ds \\ 
& \leq 
\int_0^t 2e^{-2\rho (t-s)} Q_s (Q_{t-s}( \Gamma^W(f) ))ds \\
& =
\frac{1-e^{-2\rho t}}{\rho} Q_t(\Gamma^W(f)) .
\end{align*}
This corresponds to the expected result of Item $(iii)$.

Last we prove that $(iii)$ implies $(i)$. We may use the following expansions left to the reader:
$$
Q_tf=f+tLf+ \frac{t^2}{2}L(Lf) +o(t^2)
$$
from which we deduce that
$$
Q_t (f^2) - (Q_t f)^2 = 2t \Gamma(f) + t^2[L(\Gamma(f))+2\Gamma(f,Lf)] + o(t^2).
$$
On the other hand, 
$$
2 \int_0^t Q_s \left( W^2 (Q_{t-s} f)^2\right) ds = 2tW^2f^2 + t^2[L(W^2f^2)+2W^2fLf] + o(t^2)
$$
and
$$
\frac{1-e^{-2\rho t}}{\rho} Q_t (\Gamma^W(f)) = 2t\Gamma^W(f)+t^2[-2\rho \Gamma^W(f)+2L(\Gamma^W(f))] + o(t^2).
$$
Plugging these expansions into $(iii)$ leads precisely to $(i)$. This ends the proof.
\end{proof}

In the next result, we give a condition for Inequality $(i)$ of Theorem \ref{th:g2w} to hold. Observe first that
\begin{align*}
\Gamma_2^W(f) 
& = \Gamma_2(f) + \frac{1}{2}f^2L(W^2) + W^2 \Gamma(f) + 2 W \nabla W \nabla (f^2) \\
& = \Gamma_2(f)  +
f^2 [W \Delta W + |\nabla W|^2 - W \nabla W \cdot \nabla U] + W^2 |\nabla f|^2 + 4fW \nabla W \cdot \nabla f .
\end{align*}
Set $\partial_{ij}^2$ for the second order derivative with respect to the variables $x_i$ and $x_j$.
Since $\Gamma_2(f) =\sum_{i,j=1}^n (\partial_{ij}^2 f)^2 + (\nabla f)^T (\mathrm{Hess} U) (\nabla f)$, we observe that the condition $\Gamma_2(f) \geq \rho \Gamma(f)$
is satisfied as soon as $\mathrm{Hess}(U)\geq \rho$.

\begin{theorem} \label{prop:g2}
Assume that $\Gamma_2 \geq \rho \Gamma$ for some $\rho \in \mathbb{R}$ and that
$$
\gamma := \inf_{x \in \mathbb{R}^n : W(x)\neq 0} \left(\frac{\Delta W}{W} - 3\frac{|\nabla W|^2}{W^2} - \frac{\nabla U \cdot \nabla W}{W} \right) > -\infty.
$$ 
Then, we have
\[
\Gamma^W_2(f) \geq min(\rho,\gamma)\Gamma^W(f)\]
for all $f$ smooth enough.
\end{theorem}

\noindent In the above, by convention we set $\inf \emptyset = + \infty$.

\begin{example}
Consider $U(x)=c+\frac{(1+|x|^2)^{p/2}}{p}$ and $W(x)=\frac{(1+|x|^2)^{q/2}}{q}$, $x \in \mathbb{R}^n$, $p,q \geq 1$
with $c$ so that $\int e^{-U(x)}dx=1$.
Here, as usual, $|x|=(\sum x_i^2)^{1/2}$ is the Euclidean norm. 
The (spurious) form of $U$ and $W$ is here to guarantee smoothness (indeed it would have been easier to work with $W(x)=|x|^q$ that is, however, not smooth on the whole Euclidean space).

We observe that
$\nabla U(x)=x(1+|x|^2)^{(p-2)/2}$, $\nabla W(x)=x(1+|x|^2)^{(q-2)/2}$ and $\Delta W(x)=n(1+|x|^2)^{(q-2)/2} +(q-2)|x|^2(1+|x|^2)^{(q-4)/2}$. Therefore
$$
\frac{\Delta W}{W} - 3\frac{|\nabla W|^2}{W^2} - \frac{\nabla U \cdot \nabla W}{W}
= \frac{qn}{1+|x|^2} -\frac{q|x|^2}{1+|x|^2} \left( \frac{2(q+1)}{1+|x|^2} + (1+|x|^2)^{(p-2)/2} \right)
$$
is bounded below if and only if $p \leq 2$, in which case, Theorem \ref{prop:g2} applies and leads to a non-trivial statement.

%
\end{example}

\begin{proof}[Proof of Theorem \ref{prop:g2}]
Form the expression of $\Gamma^W_2$ above, we infer that
\begin{align*}
\Gamma_2^W(f) 
& \geq  \rho |\nabla f|^2  +
f^2 [W \Delta W + |\nabla W|^2 - W  \nabla U \cdot \nabla W ] + W^2 |\nabla f|^2  + 4fW \nabla W \cdot \nabla f .
\end{align*}
Now $4fW \nabla W \cdot \nabla f \geq - 4 f^2 |\nabla W|^2 - W^2 |\nabla f|^2$ so that 
\begin{align*}
\Gamma_2^W(f) 
& \geq  \rho |\nabla f|^2  +
W^2 f^2 \left(\frac{\Delta W}{W} - 3\frac{|\nabla W|^2}{W^2} - \frac{\nabla U \cdot \nabla W}{W} \right) .
\end{align*}
The expected result follows.
\end{proof}

As an immediate corollary, we get the following useful result.

\begin{corollary} \label{cor:g2}
Assume that $\Gamma_2 \geq \rho \Gamma$ for some $\rho \in \mathbb{R}$ and that
$$
\gamma := \inf_{x \in \mathbb{R}^n : W(x)\neq 0} \left(\frac{\Delta W}{W} - 3\frac{|\nabla W|^2}{W^2} - \frac{\nabla U \cdot \nabla W}{W} \right) > - \infty .
$$ 
Then, for all $f$ smooth enough, it holds
\begin{equation} \label{eq:g2th}
\Gamma^W (Q_tf) \leq e^{-2\min(\rho,\gamma)t} Q_t( \Gamma^W(f)), \qquad t \geq 0 .
\end{equation}
\end{corollary}

\bigskip

Next, we show that the quantity $\min(\rho,\gamma)$, that appears in Theorem \ref{prop:g2} and Corollary \ref{cor:g2}, is optimal, in the sense that, for some examples of $U$ and $W$, it cannot be improved.

Observe first that, if $W \equiv 0$, then $\gamma=\infty$ and therefore  $\Gamma^W_2 \geq min(\rho,\gamma)\Gamma^W$ is equivalent to 
$\Gamma_2 \geq \rho \Gamma$ which is known to be optimal (for example for the Gaussian potential $U(x)=|x|^2/2$ for which  $\rho=1$).

In fact, consider the Gaussian potential $U(x)=\frac{|x|^2}{2} - \frac{n}{2}\log(2\pi)$ in $\mathbb{R}^n$,
and $W(x)=\sqrt{1+|x|^2}$, $x \in \mathbb{R}^n$. One has
$$
\frac{\Delta W}{W} - 3\frac{|\nabla W|^2}{W^2} - \frac{\nabla U \cdot \nabla W}{W}
= \frac{n-2}{1+|x|^2} + \frac{3}{(1+|x|^2)^2} - 1
$$
form which one infers that $\gamma=-13/12$ if $n=1$ and $\gamma=-1$ when $n \geq 2$. Since in that specific case $\rho=1$, $\min(\rho,\gamma)=\gamma$ and therefore, Corollary \ref{cor:g2} asserts that, in dimension 2 or higher,  
$\Gamma^W( Q_tf) \leq e^{2 t} Q_t(\Gamma^W(f))$, for $t \geq 0$. 
We stress that this goes in the opposite direction of \eqref{eq:be}, which, in the Gaussian setting, can be recast as
$|\nabla Q_tf|^2 \leq e^{-t} Q_t(|\nabla f|^2)$ with optimal decay $e^{-t}$. As we may prove now, 
$e^{2t}$ is also optimal.

For that purpose, consider the following family of functions
$$
f_a(x):=e^{a \cdot x}, \qquad a=(a_1,\dots,a_n),x=(x_1,\dots,x_n) \in \mathbb{R}^n 
$$ 
where as usual $a\cdot x := \sum_i a_i x_i$ is the scalar product in $\mathbb{R}^n$. Optimality can be obtained equivalently (thanks to Theorem \ref{th:g2w}) either from the bound $\Gamma^W( Q_tf) \leq e^{2 t} Q_t(\Gamma^W(f))$ (using Melher's representation formula for the Ornstein-Uhlenbeck semi-group) 
or in $\Gamma_2^W (f) \geq - \Gamma^W(f)$. We will dig on the latter by computing $\Gamma_2^W (f_a)$ and $\Gamma^W (f_a)$ for all $a,x \in \mathbb{R}^n$.

On the one hand we have
$$
\Gamma^W (f_a) = |\nabla f_a|^2 + W^2 f_a^2 = \left( |a|^2 + 1+|x|^2 \right) f_a^2 . 
$$
On the other hand,
\begin{align*}
\Gamma_2^W (f_a) 
& =
\sum_{i,j=1}^n (\partial_{ij}^2 f_a)^2 + (\nabla f_a)^T (\mathrm{Hess} U) (\nabla f_a)  +
f_a^2 [W \Delta W + |\nabla W|^2 - W \nabla W \cdot \nabla U] \\
& \quad + W^2 |\nabla f_a|^2 + 4f_aW \nabla W \cdot \nabla f_a \\
&=
f_a^2 \left( |a|^4 + 2|a|^2 + n + |x|^2(-1+|a|^2) + 4 x \cdot a \right)
\end{align*}
Therefore,
$$
\lim_{|x| \to \infty} \frac{\Gamma_2^W (f_a)}{\Gamma^W (f_a)} = -1 + |a|^2 .
$$
Finally, in the limit $|a| \to 0$, we conclude that the biggest constant $\kappa$ satisfying $\Gamma_2^W (f) \geq \kappa \Gamma^W(f)$ for all $f$ must satisfy
$\kappa \leq -1$ and therefore, by Theorem \ref{prop:g2}, $\kappa=-1$ is optimal, as announced.

\bigskip

In some situations it might be useful to deal with $\sqrt{\Gamma}$ instead of $\Gamma$.
Unfortunately, there is not a clean commutation result, as in the usual Bakry-Emery theory,
for $\Gamma^W$. However, we may prove the following proposition, that is already useful for applications. In particular, such a result was used by the authors to deal with hypercontractivity properties for some class of inhomogeneous Markov semi-groups \cite{RZ}.

\begin{proposition} \label{prop:g2sqrt}
Assume the following:\\

\noindent $(i)$ there exists $\rho \in \mathbb{R}$ such that
for all $f$ smooth enough it holds $\Gamma_2(f) \geq \rho \Gamma(f)$;\\

\noindent $(ii)$ $c:=\max\left( 2 \| |\nabla W| \|_\infty , \sup_{x: W(x) \neq 0} \left( \frac{LW}{W} - \rho \right)_- \right) < \infty$. \\
Then, for all $f$ non-negative, it holds
$$
\sqrt{\Gamma(P_tf)} + W P_tf 
\leq e^{(c-\rho)t}  P_t \left( \sqrt{\Gamma(f)} + W f \right).
$$
\end{proposition}

\begin{proof}
Following Bakry-Emery, see \cite[proof of Proposition 5.4.5]{ane}, introduce
$\Psi(s)=e^{-\rho s} P_s \left( \sqrt{\Gamma(P_{t-s}f)} + WP_{t-s}f \right)$, $s \in [0,t]$, $t$ being fixed. Therefore, setting $g:=P_{t-s}f$, one has
\begin{align*}
\Psi'(s)= - \rho \Psi(s)+e^{-\rho s} P_s L \sqrt{\Gamma(g)} + P_s(L(Wg))
-e^{-\rho s}P_s \left( \frac{\Gamma(g,Lg)}{\sqrt{\Gamma(g)}}\right)
-e^{-\rho s}P_s \left( W Lg \right) .
\end{align*}
Now
$$
L \sqrt{\Gamma(g)} = \frac{L \Gamma(g)}{2 \sqrt{\Gamma(g)}} - \frac{\Gamma(\Gamma(g))}{4 \Gamma(g)^{3/2}}
$$ 
Hence, after some algebra, we get
$$
\Psi'(s)=
\frac{e^{-\rho s}}{4} P_s \left( \frac{4\Gamma(g)(\Gamma_2(g)-\rho \Gamma(g))- \Gamma(\Gamma(g))}{\Gamma(g)^{3/2}} \right)
+ e^{-\rho s} P_s \left( L(Wg) - W Lg - \rho Wg \right) .
$$
Assumption $(i)$ ensures that the first term of the right hand side of the latter is non-negative (see \cite[Lemma 5.4.4]{ane}). On the other hand,
$$
L(Wg) - W Lg - \rho Wg 
= 
g(LW - \rho W) + 2 \nabla W \cdot \nabla g
\geq 
-c \left( |\nabla g| + Wg \right) .
$$
It follows that
$\Psi'(s) \geq -c \Psi(s)$. In turn, $\Psi(t) \geq \Psi(0)e^{-ct}$ from which the desired result follows.
\end{proof}

\bibliographystyle{alpha}
\bibliography{gamma2-RZ}

\newcommand{\etalchar}[1]{$^{#1}$}
\def\cprime{$'$}
\begin{thebibliography}{ABC{\etalchar{+}}00}

\bibitem[ABC{\etalchar{+}}00]{ane}
C.~An\'e, S.~Blach\`ere, D.~Chafai, P.~Foug\`eres, I.~Gentil, F.~Malrieu,
  C.~Roberto, and G.~Scheffer.
\newblock {\em Sur les in\'egalit\'es de {S}obolev logarithmiques.}, volume~10
  of {\em Panoramas et {S}ynth\`eses}.
\newblock S.M.F., Paris, 2000.

\bibitem[Bak94]{bakry}
D.~Bakry.
\newblock L'hypercontractivit\'e et son utilisation en th\'eorie des
  semigroupes.
\newblock In {\em Lectures on probability theory. \'Ecole d'été de
  probabilités de St-Flour 1992}, volume 1581 of {\em Lecture Notes in Math.},
  pages 1--114. Springer, Berlin, 1994.

\bibitem[BE85]{bakry-emery}
D.~Bakry and M.~Emery.
\newblock Diffusions hypercontractives.
\newblock In {\em S\'eminaire de probabilit\'es, XIX, 1983/84}, pages 177--206.
  Springer, Berlin, 1985.

\bibitem[BGL14]{BGL}
D.~Bakry, I.~Gentil, and M.~Ledoux.
\newblock {\em Analysis and geometry of {M}arkov diffusion operators}, volume
  348 of {\em Grundlehren der Mathematischen Wissenschaften [Fundamental
  Principles of Mathematical Sciences]}.
\newblock Springer, Cham, 2014.

\bibitem[RZ]{RZ}
C.~Roberto and B.~Zegarlinski.
\newblock Hypercontractivity for markov semi-groups.
\newblock preprint arXiv:2101.01616.

\end{thebibliography}

\end{document}